\documentclass[12pt]{article}
\usepackage{amsxtra,amssymb,amsthm,amsmath,latexsym}

\theoremstyle{plain}

\def\R{{\mathbb R}}

\def\oH{\buildrel\circ\over H}
\def\oH1{\buildrel\circ\over H\kern-.02in{}^1}
\def\qed{{\hfill $\Box$}}

\begin{document}

\title{ The shape of the ear canal
   \thanks{key words:  acoustic waves, 
inverse problems, ear
    }
}

\author{
A.G. Ramm\\
Mathematics Department, 
Kansas State University, \\
 Manhattan, KS 66506-2602, USA\\
ramm@math.ksu.edu\\
}

\date{}

\maketitle\thispagestyle{empty}

\begin{abstract} 
It is proved that the measurement of the acoustic pressure 
on the ear membrane allows one to determine the shape of the ear canal 
uniquely.
\end{abstract}


\section{Introduction}
Consider a bounded domain $D \subset \R^n$, $n = 3,$ with a
Lipschitz  boundary $S$. Let $F$ be an open subset on $S$, a membrane, 
$G=S\setminus F$,  $\Gamma=\partial F$, and $N$ is the outer unit normal 
to $S$. The domain $D$ models the ear canal. The acoustic pressure $u$
is given on the membrane $F$, $u=f\not\equiv 0$ on $F$.  This pressure is 
zero on $G$, and the 
normal component of the velocity, $u_N:=h$, is measured on $F$.
The problem we are concerned with is: {\it given the datum $\{f,h\}$,
with $f\not\equiv 0$, for a single $f$, can one recover uniquely
the shape of $D$, the ear canal?}

Let us introduce now the 
corresponding mathematical formulation of the above inverse problem. 
Consider the problem:
$$\nabla^2 u+k^2 u=0 \hbox{\ in\ } D, \quad
  u = f \hbox{\ on\ } F,  \quad u = 0 \hbox{\ on\ } G. \eqno{(1.1)}$$
We assume that $k^2$ is not a Dirichlet eigenvalue of the
Laplacian in $D$. This assumption will be removed later.
If this assumption holds,  then the solution to problem
(1.1) is unique. Thus, its normal derivative, $u_N:=h$ on $F$, is uniquely 
determined. Suppose one can measure $h$ on $F$ for a single $f\in C^1(F), 
\,\,f\not\equiv 0$.

The inverse problem (IP) we are interested in can now be formulated:
 
{\it Does this datum determine $G$ uniquely?}

Thus, we assume that $F$, $f$ and $h$ are known, that $k^2$ is not a 
Dirichlet eigenvalue of the Laplacian in $D$,  and we want 
to determine the unknown part $G$ of the boundary $S$. 

Let $\Lambda$ be
the smallest eigenvalue of the Dirichlet Laplacian $L$ in $D$.
Let us assume that 
$$\Lambda>k^2.
\eqno{(1.2)}$$
Then problem (1.1) is uniquely solvable.
Assumption (1.2) in our problem is practically not a serious restriction, 
because the wavelength in our experiment can be chosen as we wish.
Since the upper bound on the width $d$ of the ear canal is known, 
and since 
$$\Lambda>\frac 1 { d^2}, \eqno{ (1.3)}$$ 
one can choose $k^2<\frac 
{1} {d^2}$ to satisfy assumption (1.2). A proof of the
estimate  (1.3) is given at the end of this note.

We discuss the Dirichlet condition but 
a similar argument is applicable to the Neumann and Robin boundary 
conditions. Boundary-value problems and scattering problems in rough 
domains were studied in [1].

Our basic result is the following theorem:
  
{\bf Theorem 1.} {\it If (1.2) holds then the above data determine $G$ 
uniquely}.

{\bf Remark 1.}  If $k^2$ is an eigenvalue of the Dirichlet
Laplacian $L$ in $D$, and $m(k)$ is the total multiplicity 
of the spectrum of $L$ on the semiaxis $\lambda\leq k^2$,
then $G$ is uniquely defined by the data $\{f_j, h_j\}_{1\leq j \leq 
m(k)+1}$, where $\{f_j\}_{1\leq j \leq m(k)+1}$ is an arbitrary fixed 
linearly independent system of functions in $C(F)$.

In Section 2 proofs are given, in Section 3 a numerical approach to 
computing $G$ is discussed very briefly, and in Section 4 conclusions are 
formulated.

\section{Proofs.}

{\bf Proof of Theorem 1.}   

Suppose that there are two surfaces $G_1$ and $G_2$, which generate the 
same data, that is, given $f$ on $F$ they generate the same function $h$ 
on $F$.  Let $D_1, u_1$ and 
$D_2, u_2$ be the corresponding domains and 
solutions to (1.1). Denote $w:=u_1-u_2$, $D^{12}:=D_1\cap D_1$, 
$D_{12}:=D_1\cup D_2$, $D_3:=D_1\setminus D^{12}$, $D_4:=D_2\setminus 
D^{12}$. Note that $ w= w_N = 0$ on $F$, since the data 
$f$ and $h$ are the same by our assumption.
Consequently, one has:
$$\nabla^2 w+k^2w=0 \hbox{\ in\ } D^{12}, \quad
  w= w_N = 0 \hbox{\ on\ } F. \eqno{(2.1)}$$
By the uniqueness of the solution to the Cauchy problem for elliptic 
equations, one concludes that $w=0$ in $ D^{12}$. Thus, $u_1=u_2=0$
on $\partial D^{12}$, and $u_1=0$ on $\partial D_3$. Therefore
$$\nabla^2 u_1+k^2u_1=0 \hbox{\ in\ } D_3, \quad
  u_1 = 0 \hbox{\ on\ } \partial D_3. \eqno{(2.2)}$$
Since $D_3\subset D$, it follows that $\Lambda (D_3)>\Lambda (D)>k^2$.
Therefore $k^2$ is not a Dirichlet eigenvalue of the Laplacian in 
$D_3$, so $u_1=0$ in $D_3$, and, by the unique continuation property
for solutions of the homogeneous Helmholtz equations one concludes that 
$u_1=0$ in $D_1$. In particular, $u_1=0$ on $F$, which is a 
contradiction, since $u_1=f\neq 0$ on $F$ by the assumption. 
Theorem 1 is proved. \qed

{\bf Proof of Remark 1.} Suppose that $k^2>0$ is arbitrarily fixed, and 
the data are $\{f_j, h_j\}_{1\leq j\leq m(k)}$.
Using the same argument as in the proof of Theorem 1, one arrives at the 
conclusion (2.2) with $u_{j,1}$ in place of $u_1$, where $u_{j,1}$
solves (1.1) with $f=f_j$, $1\leq j \leq m(k)+1$.
Since the total multiplicity of the spectrum of the Dirichlet Laplacian 
in $D$ is not more that $m(k)$,  one  concludes that 
$D_1=D_2$. Remark 1 is proved. \qed


{\it Proof of estimate (1.3)}. 

For convenience of the reader and to make the presentation self-contained, 
we give a short proof 
of inequality (1.3), although this inequality is known. 
Let $\alpha$ be a unit vector, and
$d(\alpha)$ be the width of $D$ in the direction $\alpha$, that is, the 
distance between two planes, tangent to the boundary $S$ of $D$
and perpendicular to the vector $\alpha$, so that $D$ lies between these 
two planes. Let 
$$d:=\min_{\alpha}d(\alpha)>0.$$
By the variational definition of $\Lambda$ one has:
$\Lambda=\min \int_D|\nabla u|^2 dx,$ where the minimization is
taken over all $u\in H^1$, vanishing on $S$ and normalized, 
$||u||_{L^2(D)}=1$. Denote $s:=x_1, \, y:=(x_2,x_3),$
and choose the direction of $x_1-$axis along the direction $\alpha$,
which minimizes $d(\alpha)$, so that the width of $D$ in the direction 
of axis $x_1$ equals $d$. Then one has:
$$u(s,y)=\int_a^s u_t(t,y)dt,$$ so 
$$|u(s,y)|^2\leq \int_a^s |u_t(t,y)|^2dt 
(s-a)\leq d \int_a^b |u_t(t,y)|^2dt,$$ 
where $s=a$ and $s=b$ are the 
equations of the two tangent to $S$ planes, the distance between them is 
$d=b-a$, and $D$ is located between these planes.

Denote by $F_s$ the crossection of $D$ by the plane $x_1=s$, $a<s<b$.
Integrating the last inequality with respect to $y$ over $F_s$, 
and then with respect to $s$ between $a$ and $b$, one gets:
$$||u||^2_{L^2(D)}\leq d^2 ||\nabla u||^2,$$  
  which implies inequality (1.3). \qed

\section{ A discussion of the numerical aspects of the 
problem}
Let us discuss very briefly how to calculate $G$ numerically, 
given the data $\{f,h\}$. We have to calculate three scalar functions
which give a parametric equation of $G$, of the form $x_j=\phi_j(t,v)$,
$j=1,2,3$, $t,v\in (0,1)\times (0,1)$ are parameters, the three functions
$ \phi_j$ are unknown. They are to be found from the three conditions:
$u=f$ on $F$, $u_N=h$ on $F$, and $u=0$ on $G$. One writes a representation 
for $u$ in $D$ by Green's formula
$$u=\int_F[g(x,s)h(s)-g_N(x,s)f(s)]ds+\int_Gg(x,s)u_N(s)ds:=U(x)+V(x),$$
where $g=\frac {e^{ik|x-s|}}{4\pi|x-s|}$, the function $u_N$ on $G$ is not 
known, $ds$ is an element of the surface area, the functions $\phi_j$,
$j=1,2,3,$ are not known, the function $U(x)$ is known, and $V(x)$ is a 
single-layer potential with the unknown density $u_N:=H$ on $G$. To 
determine 
the four unknown functions $\phi_j, j=1,2,3,$ and $H$, one has three 
boundary equations:
$u=f$ on $F$, $u_N=h$ on $F$, $u=0$ on $G$, and 
the fourth equation for $H$ on $G$ one obtains from the formula 
$H=U_N+\frac 1 2 AH+\frac 1 2 H$ for the normal derivative of the 
single-layer 
potential on $G$, which yields an integral equation for $H$:
$H=AH+U_N$. Here $AH:=2\int_G g_N(s,s')H(s')ds'$.
One expresses $ds$ in terms of the functions $\phi_j$, and get
four nonlinear integral equations for four unknown functions
$\phi_j, \,j=1,2,3,$ and $H$.
These equations do have a solution if the data $\{f,h\}$ are exact.
The solution can be obtained by a Newton-type method. The problem is 
ill-posed, because small variations of $h$ lead to a function which 
may be not a normal derivative on $F$ of the solution to problem (1.1) 
corresponding to the given function $f$. Thus, one has to use
a regularization for solving the above integral equations numerically.

\section{Conclusion}

The basic result of this note is the proof of the following statement:
if one applies some pressure to the eardrum and measures the corresponding 
normal component of the velocity of this drum, then one can uniquely 
determine from these data the shape of the ear canal.
This note is of theoretical nature, but hopefully it may lead to a
progress in the construction of better hearing devices.

\end{document}